\documentclass[12pt]{article}
\usepackage{amsmath}
\usepackage{amssymb}
\usepackage{url}

\DeclareMathOperator*{\newand}{{\&}}
\date{}

\begin{document}
The original article (in Russian) appeared in:
\emph{Teoriya Algorifmov i Matematicheskaya Logika }(a collection of papers
dedicated to
A.\,A.\,Mar\-kov), Vychislitel'nyi Tsentr Akademii Nauk SSSR, Moscow,
1974, pages 112--123;
\url{http://www.ams.org/mathscinet-getitem?mr=0406780}.

The original article was reprinted in
\emph{The Collected Works of Julia Robinson},
\url{http://www.ams.org/mathscinet-getitem?mr=1411448}.

\vspace{20mm}

\hfill\begin{tabular}{cc}
{\sc Yuri Matiyasevich\thanks{
The English translation was carried out
by M.\,A.\,Vsemirnov and
 edited by J.\,P.\,Jones and M.\,Davis; the first author is
 grateful to all three of them.},}&{\sc Julia Robinson}\\
{\it Leningrad}&{\it Berkeley, USA}
\end{tabular}
\vspace{5mm}
\begin{center}
TWO UNIVERSAL 3-QUANTIFIER REPRESENTATIONS OF RECURSIVELY
ENUMERABLE SETS
\end{center}

\vspace{15mm}

{\bf 1.} We shall use the following notation: Lower-case Latin
letters from $a$ to $n$ (inclusive) with or without subscripts
will be used as variables for nonnegative integers, the remaining lower-case
Latin letters will be used as variables for integers. Analogously,
lower case Greek letters from $\alpha$ to $\nu$ will be used as
metavariables for nonnegative integers, and the remaining the Greek letters
will be used as
metavariables for integers.

Upper case Latin letters will denote polynomials.
Here and below it is to be understood that only polynomials with integer coefficients are being considered.

{\bf 2.} We say that a set $\mathfrak R$  of nonnegative integers
is represented by an arithmetic formula $\mathfrak F$ with one
free variable $a$ if the  equivalence $a\in \mathfrak R\Leftrightarrow \mathfrak F$
is true.

As K.\,G{\"o}del showed, any recursively enumerable set
is represented by an arithmetic formula. One can improve this
result by restricting the kinds of formulas in various ways.
Such restricted representations were investigated in [2-13]. The aim of this paper is to show  that every recursively enumerable set is represented by formulas
of each of the two kinds following:
\begin{gather}
\exists b \exists c \newand_{\iota=1}^{\varepsilon}
\exists d [
P_{\iota}(a,b,c)<D_{\iota}(a,b,c)d<Q_{\iota}(a,b,c) ],
%\eqno (1)
\\
\exists b \exists c \forall f
[f\le F(a,b,c) \Rightarrow W(a,b,c,f)>0].
%\eqno (2)
\end{gather}

{\bf 3.} Let $\mathfrak R$ be a recursively enumerable set of
non-negative integers. We begin with a formula that represents the set $\mathfrak R$
of the form
\begin{equation}
\exists h_1 \dots  \exists h_{\delta} [R(a,h_1,\dots, h_{\delta})=0],
%\eqno (3)
\end{equation}
 (the existence of such a formula
is proved, for example, in [6-9,14]).

Denoting the degree of the polynomial $R$ by $\lambda$,
without loss of generality we may assume that $\lambda\ge 1$.

In order that formula (2) be equivalent to formula (3), the pair
$\left\langle b,c \right\rangle$,
whose existence is asserted in (2), must carry
all of the  information  contained in the $\delta$\/-tuple
$\left\langle h_1,\dots h_{\delta} \right\rangle$,
whose existence is asserted in~(3).

Many methods are known for coding tuples of nonnegative integers
using a single nonnegative integer or a pair of such integers.
The rather unusual method that we use allows us
to check the truth of the relation
\begin{equation}
R(a,h_1,\dots, h_{\delta} )=0
%\eqno (4)
\end{equation}
directly from the code, without first finding the
individual numbers $h_1,\dots, h_{\delta}$.
We define  $B(h_1,\dots,h_{\delta},k)$ to be the polynomial
$$
\sum_{\iota=1}^{\delta} h_\iota k^{(\lambda+1)^{\iota}}.
$$
This polynomial has the ``geometric'' interpretation:
if $k$ is greater than each of the numbers $1, h_1,\dots,h_{\delta}$,
then
$h_1,\dots,h_{\delta}$ are the corresponding
$(\lambda+1)$-th,\dots, $(\lambda+1)^{\delta}$-th
 digits of the number
$B(h_1,\dots,h_{\delta},k)$
in the $k$\/-ary  number system, while
all of the other digits are zeros.

One can easily verify that an identity of the following type holds:
\begin{multline}
(1+ak+B(h_1,\dots, h_{\delta},k))^\lambda=\\
\sum_{\alpha_0+\dots+\alpha_\delta\le \lambda}
\kappa_{\alpha_0,\dots,\alpha_\delta}
a^{\alpha_0} h_1^{\alpha_1}\dots h_\delta^{\alpha_\delta}
k^{N(\alpha_0,\dots,\alpha_\delta)}
%%%%%                      ^^^^^^^
%%%%%                     \iota or 1 in the original
%%%%%%
%\eqno (5)
\end{multline}
where
$$
N(l_0,\dots, l_\delta)=\sum_{\iota=0}^\delta l_\iota (\lambda+1)^\iota,
$$
and
$\kappa_{\alpha_0,\dots,\alpha_\delta}$ are positive integers.
One can easily see that this polynomial $N$ has the following property:
\begin{multline}
\biggl\{
\newand_{\iota=0}^\delta (( l^{'}_\iota\le \lambda)\, \& \,
(l^{''}_\iota\le \lambda) \, \& \,
%\\
 N(l^{'}_0,\dots, l_\delta^{'}) =
N(l^{''}_0,\dots, l_\delta^{''}) \biggr\}
\Rightarrow
\newand_{\iota=0}^\delta (l_\iota^{'}=l_\iota^{''}).
%\eqno (6)
\end{multline}
This property is obvious applying the above mentioned ``geometric'' interpretation
to the polynomial $N$: if $l_\iota\le \lambda$, then
$l_\delta,\dots, l_0$ are precisely the digits  in the expansion
of $N(l_0,\dots, l_\delta)$
in the number system with base $\lambda+1$.
The polynomial $B$ was chosen in such a way that in the $k$\/-ary
expansion
of the number
$B(h_1,\dots,h_{\delta},k)$ the non-zero digits  are placed
in special locations in order
to obtain property~(6).
Property~(6) allows us to give ``a geometric interpretation''
of identity~(5):
in the $k$\/-ary expansion of the number
$(1+ak+B(h_1,\dots, h_{\delta},k))^\lambda$
the digits are all possible numbers of the form
\begin{equation}
\kappa_{\alpha_0,\dots,\alpha_\delta}
a^{\alpha_0} h_1^{\alpha_1}\dots h_\delta^{\alpha_\delta}
%\eqno (7)
\end{equation}
provided that $k$ exceeds each of them.

Without loss of generality we shall assume that the polynomial $R$
is a linear combination of monomials~(7):
\begin{equation}
R(a,h_1,\dots, h_{\delta}))=
\sum_{\alpha_0+\dots+\alpha_\delta\le \lambda}
\rho_{\alpha_0,\dots,\alpha_\delta}
\kappa_{\alpha_0,\dots,\alpha_\delta}
a^{\alpha_0} h_1^{\alpha_1}\dots h_\delta^{\alpha_\delta}
%\eqno(8)
\end{equation}

Obviously, if $l_1+\dots+l_\delta\le \lambda$, then
$$
N(l_1,\dots,l_\delta)\le \lambda(\lambda+1)^\delta.
$$
Let us denote $\lambda(\lambda+1)^\delta$ by $\nu$, and
the polynomial
$$
\sum_{\alpha_0+\dots+\alpha_\delta\le \lambda}
\rho_{\alpha_0,\dots,\alpha_\delta}
k^{\nu-N(\alpha_0,\dots,\alpha_\delta)}
$$
by $V(k)$.
One can easily see that an identity of the following type holds:
\begin{equation}
V(k) (1+ak+B(h_1,\dots, h_{\delta},k))^\lambda
=\sum_{\iota=0}^{2\nu}
T_\iota(a,h_1,\dots, h_{\delta}) k^\iota,
%\eqno (9)
\end{equation}
where $T_0$,\dots,$T_{2\nu}$ are polynomials whose degrees
do not exceed $\lambda$.

One can interpret identity (9) in a natural manner
if one considers a $k$\/-ary number system in which negative
digits are allowed; for instance, one may require that
\begin{equation}
k>|2T_\iota(a,h_1,\dots, h_{\delta})|, \qquad (\iota=0,\dots,2\nu)
%%%%       ^^^
%%%%       a_1 in the original
%%%%
%\eqno(10)
\end{equation}
and consider the system with digits ranging from
$[-(k-1)/2]$ to \mbox{$[(k-1)/2]$}.

It is easy to check  that (5), (6), (8) and (9) imply the identity
\begin{equation}
T_\nu(a,h_1,\dots, h_{\delta})=
R(a,h_1,\dots, h_{\delta}).
%\eqno(11)
\end{equation}
Thus, if
\begin{equation}
b=B(h_1,\dots, h_{\delta},k)
%\eqno(12)
\end{equation}
and $k$ is sufficiently large that the inequalities~(10) are satisfied, then
the relation~(4)  holds if and only if the digit in the $\nu$\/-th place
in the $k$\/-ary expansion of the number $V(k)(1+ak+b)^\lambda$ is zero.
As we shall show below, the latter condition can be easily written
using a single existential quantifier.

\

\noindent{\bf Lemma 1.}
For any $a$, $b$, $h_1,\dots,h_\delta$, $k$ satisfying conditions~(10) and~(12),
the relation (4) holds if and only if there exists an integer $z$ such that
\begin{equation}
-k^\nu< 2(V(k)(1+ak+b)^\lambda-z k^{\nu+1}) < k^\nu.
%%%%                                ^^^^
%%%%                                \lambda in the original
%%%%
%\eqno (13)
\end{equation}

\

\noindent{\bf Necessity.}
Put
$$
z=\sum_{\iota=\nu+1}^{2\nu} T_\iota(a,h_1,\dots, h_{\delta}) k^{\iota-\nu-1}
$$
By (9), (12) and (4),
$$
V(k)(1+ak+b)^\lambda-z k^{\nu+1} =
\sum_{\iota=0}^{\nu-1} T_\iota(a,h_1,\dots, h_{\delta}) k^{\iota}.
$$
We deduce from (10) that
\begin{multline}
\left|
2\sum_{\iota=0}^{\nu-1} T_\iota(a,h_1,\dots, h_{\delta}) k^{\iota}
\right| \le \\
\le \sum_{\iota=0}^{\nu-1} |2T_\iota(a,h_1,\dots, h_{\delta}) k^{\iota}|
 \le \\
\le \sum_{\iota=0}^{\nu-1} (k-1) k^{\iota} =k^\nu-1<k^\nu,
%\eqno (14)
\end{multline}
so that inequalities (13) are satisfied.

\noindent{\bf Sufficiency.}
It is easy to see that there exists at most one integer $y$ such that
$$
-k^\nu< 2(V(k)(1+ak+b)^\lambda-y k^{\nu}) < k^\nu.
$$
On  the one hand, by (13), $y$ equals  $zk$, while on the other  hand
$y$ equals
$$
\sum_{\iota=\nu}^{2\nu} T_\iota(a,h_1,\dots, h_{\delta}) k^{\iota-\nu}
$$
since , by (9) and (12),
\begin{multline*}
V(k)(1+ak+b)^\lambda-
\biggl(\sum_{\iota=\nu}^{2\nu} T_\iota(a,h_1,\dots, h_{\delta}) k^{\iota-\nu}
\biggr)k^\nu=\\
=\sum_{\iota=0}^{\nu-1} T_\iota(a,h_1,\dots, h_{\delta}) k^{\iota}
\end{multline*}
and inequalities (14) hold. Thus,
$$
\sum_{\iota=\nu}^{2\nu} T_\iota(a,h_1,\dots, h_{\delta}) k^{\iota-\nu}=zk.
$$
Passing from this equation to a congruence, we get
$$
T_\nu(a,h_1,\dots, h_{\delta}) \equiv 0 \pmod k .
$$
This, together with (10), gives us the equality
\begin{equation}
T_\nu(a,h_1,\dots, h_{\delta}) = 0  .
%\eqno (15)
\end{equation}
Now (4) follows from (11) and (15).

The Lemma is proved.

{\bf 4.} We now proceed to transform the inequalities (10).
Let $\gamma$ be a positive integer that exceeds twice
the  sum of the absolute values of the coefficients of all the
polynomial $T_1,\dots, T_{2\nu}$.
Obviously, the following inequalities hold:
$$
|2T_\iota(a,h_1,\dots, h_{\delta}) | < \gamma
(\max \{1,a,h_1,\dots, h_{\delta}\})^\lambda.
$$
Thus, if
$$
c \ge \max \{h_1,\dots, h_{\delta}\}
$$
and
$$
k=K(a,c)=\gamma(2+a+c)^\lambda,
$$
then inequalities (10) are satisfied.

Using Lemma 1 one can easily show that formula (3) is equivalent to the
formula
$$
\exists b \exists c [{\mathfrak F}_1\, \&\, \exists z {\mathfrak F}_2],
$$
where, here and below ${\mathfrak F}_1$ denotes  the formula
$$
%\begin{multline*}
\exists h_1
\dots
\exists h_\delta [
c \ge \max \{h_1,\dots, h_{\delta}\} \, \& \,
b=H(h_1,\dots, h_\delta, K(a,c))]
%\end{multline*}
$$
and ${\mathfrak F}_2$ denotes  the formula
\begin{multline*}
-(K(a,c))^\nu< 2(V(K(a,c))(1+aK(a,c)+b)^\lambda-\\
-       z K(a,c)^{\nu+1}) < (K(a,c))^\nu.
\end{multline*}

\

\noindent{\bf Lemma 2.} Formula $\mathfrak F_1$ is equivalent to
the formula
\begin{equation}
{\mathfrak F}_3 \,\&\,
{\mathfrak F}_4 \,\&\,
{\mathfrak F}_5 ,
%\eqno(16)
\end{equation}
where, here and below
${\mathfrak F}_3$
denotes the formula
$$
\exists d [b=d(K(a,c))^{\lambda+1}],
$$
${\mathfrak F}_4$
denotes the formula
%\begin{multline*}
$$
\newand_{\iota=1}^{\delta-1}
\exists d
\exists e
[b=d(K(a,c))^{(\lambda+1)^{\iota+1}}+e \,\& \,
e<(c+1)(K(a,c))^{(\lambda+1)^{\iota}} ],
%\end{multline*}
$$
and ${\mathfrak F}_5$
denotes the formula
$$
b<(c+1)(K(a,c))^{(\lambda+1)^{\delta}}.
$$

\

The truth of this lemma becomes quite clear, if one notes that
each of the  formulas ${\mathfrak F}_1$ and~(16) mean that  in the
expansion of the number b in the number system with the base $K(a,c)$
the non-zero digits can only occupy the $(\lambda+1)$\/-th,\dots,
$(\lambda+1)^\delta$\/-th positions, and moreover, these
digits do not exceed $c$.

{\bf 5.} Combining Lemmas 1 and 2. we see that formula (3) is equivalent to the formula
\begin{equation}
\exists b \exists c [ \exists z {\mathfrak F}_2 \, \&\,
{\mathfrak F}_3 \, \&\,
{\mathfrak F}_4 \, \&\,
{\mathfrak F}_5 ].
%\eqno(17)
\end{equation}

\

\noindent{\bf Theorem 1.} Every recursively enumerable set of nonnegative
integers can be represented by a formula of the form (1).

\

\noindent{\bf Proof.} The desired formula can be obtained
from the formula~(17) by means of easy algebraic transformations.

Formula ${\mathfrak F}_2$ contains the variable $z$, whose possible values
are all integers. However,  it follows from ${\mathfrak F}_2$ that
$$\
%begin{multline*}
2z(K(a,c))^{\nu+1}>-(K(a,c))^\nu
+ 2V(K(a.c))(1+aK(a,c)+b)^\lambda,
%\end{multline*}
$$
so that
$$
z\ge V(K(a.c))(1+aK(a,c)+b)^\lambda.
$$
Let us denote by ${\mathfrak F}_6$ the formula,
which is obtained from
${\mathfrak F}_2$  by substituting the polynomial
$$
d+ V(K(a,c))(1+aK(a,c)+b)^\lambda
$$
for $z$ and by transposing terms, which
do and do not contain $d$ to opposite sides of the inequalities.
Obviously, the formula
$\exists z {\mathfrak F}_2$
is equivalent to the formula $\exists d {\mathfrak F}_6$.

We transform the formula
${\mathfrak F}_3$ into the equivalent formula
\begin{equation}
\exists d [ b-1 < (K(a,c))^{\lambda+1} d< b+1].
%\eqno (18)
\end{equation}

Each of the conjuncts composing the formula
${\mathfrak F}_4$ includes an equation that enable us
to express
$e$ explicitly in terms of $a$, $b$, $c$ and $d$,
and, therefore, to eliminate this variable.
In addition, we must impose an inequality to insure
the non-negativity of $e$. Finally,
we obtain the formula
\begin{equation}
\newand_{\iota=1}^{\delta-1}
\exists d [ b-(c+1)(K(a,c))^{(\lambda+1)^\iota}<
(K(a,c))^{(\lambda+1)^{\iota-1}} d <b+1].
%\eqno(19)
\end{equation}
Finally we must replace the formula ${\mathfrak F}_5$ by
the equivalent formula
$$
\exists d [ b-1 <bd < (c+1)(K(a,c))^{(\lambda+1)^\delta}].
$$

The theorem is proved.

{\bf 6.} Now we turn to constructing a formula of the form~(2)
that represents the set ${\mathfrak R}$.
For this purpose we first show that for any $\epsilon$
there exist polynomials $F_\epsilon$ and $W_\epsilon$ in
$2\epsilon$ and
$3\epsilon+1$ variables respectively, such that:
if  the numbers
$g_1$,\dots,$g_\epsilon$,
$s_1$,\dots,$s_\epsilon$,
$t_1$,\dots,$t_\epsilon$ satisfy the inequalities
\begin{equation}
0<g_\iota, \quad t_\iota-s_\iota\le g_\iota \quad (\iota=1,\dots,\epsilon),
%\eqno(20)
\end{equation}
then the
formula
\begin{equation}
\newand_{\iota=1}^{\epsilon} \exists z [s_\iota <z g_\iota <t_\iota]
%\eqno(21)
\end{equation}
is equivalent to the formula
\begin{multline*}
\forall f [f\le F_\epsilon(s_1,\dots,s_\epsilon,t_1,\dots, t_\epsilon)
\Rightarrow \\
W_\epsilon(g_1,\dots,g_\epsilon,
s_1,\dots,s_\epsilon,t_1,\dots, t_\epsilon,f)>0].
\end{multline*}
We start with the case $\epsilon=1$ and find, to begin with, polynomials
$X$ and $Y$ such that for $g>0$  the formula
\begin{equation}
\exists z [s <z g <t]
%\eqno(22)
\end{equation}
is equivalent to the formula
\begin{multline}
\forall y [-s^2-t^2-2<y\le s^2+t^2+2
\Rightarrow  X(g,s,t,y)>0 \vee \\
Y(g,s,t,y) >0].
%\eqno(23)
\end{multline}
Lacking existential quantifiers, formula (23) must somehow contain
complete information about an integer $z$ that satisfies the
inequalities
\begin{equation}
s<zg<t.
%\eqno(24)
\end{equation}

We shall verify equivalence between formulas of the forms (23) and (22)
by means of the following obvious lemma, which may be regarded
as a discrete analogue of the Cauchy theorem about the vanishing
of a continuous function, whose values at the endpoints
of an interval have opposite signs.

Let $p$ and $q$ be integers such that $p<q$,
let $\Phi$ and $\Psi$ be unary predicates defined for all
integers between $p$ and $q$.
If $\Phi(p)\,\&\,\Psi(q)$ holds and for any $w$, such that $p<w<q$,
$\Phi(w)\,\vee \,\Psi(w)$ holds, then
there exists an integer $r$ such that
$p\le r\le q$ and
$\Phi(r)\,\& \,\Psi(r+1)$.

\

\noindent{\bf Lemma 3.} If
\begin{equation}
g>0,
%\eqno(25)
\end{equation}
then formula (22) is equivalent to the formula
\begin{multline}
\forall y [-s^2-t^2-2<y\le s^2+t^2+2
\Rightarrow  %\\
(y-1)g-s >0 \vee t-yg>0].
%\eqno(26)
\end{multline}

\

\noindent{\bf Proof.} Let $g$, $s$, $t$ satisfy conditions~(25) and~(26).
We will show that they satisfy condition~(22), as well.

By (25),
\begin{gather*}
t-(-s^2-t^2-1)g\ge t+s^2+t^2+1>0,\\
(s^2+t^2+1)g-s\ge s^2+t^2+1-s>0.
\end{gather*}
By the discrete analogue of the Cauchy theorem mentioned above,
we have that there exists $z$ such that
$$
t-zg>0 \ \& \ zg-s>0.
$$

Thus, condition (22) is satisfied.

Now, let $g$, $s$ and $t$ satisfy conditions (25) and (22).
We will find a $z$ that satisfies inequalities (24). Suppose that
condition (26) doesn't hold. Let $y$ be a number such that
\begin{equation}
(y-1)g-s\le 0 \,\&\, t-yg\le 0.
%\eqno(27)
\end{equation}
From (24) and (27) we obtain
$$
(y-1)g\le s< zg, \quad zg <t<yg.
$$
Consequently
$$
y-1<z<y.
$$
This contradiction  completes the proof of the equivalence
of formulas (22) and (26).

Note, that  if
\begin{equation}
t-s\le g,
%\eqno(28)
\end{equation}
then two inequalities in formula (27) are inconsistent.
Moreover, if $(y-1)g-s>0$, then $t-yg<0$, and conversely
if $t-yg>0$, then $(y-1)g-s<0$.
This enables us to transform the disjunction of a pair of inequalities
into a single one:
\begin{multline*}
(y-1)g-s>0 \vee t-yg>0 \Leftrightarrow
((y-1)g-s>0\,\&\, t-yg<0)) \,\vee        \\
 (t-yg>0\,\&\,(y-1)g<0)) \Leftrightarrow ((y-1)g-s)(yg-t)>0.
\end{multline*}

Thus, if inequalities (25) and (28) are satisfied,
then formula (22) is equivalent to the formula
$$
\forall y [-s^2-t^2-2<y\le s^2+t^2+2 \Rightarrow Z(g,s,t,y)>0],
$$
where, here and below $Z(g,s,t,y)$ denotes the polynomial
$$
((y-1)g-s)(yg-t).
$$

Note, that if $g>0$, then
\begin{equation}
\forall y [y\le -s^2-t^2-2\, \vee\, y> s^2+t^2+2 \Rightarrow Z(g,s,t,y)>0].
%\eqno(29)
\end{equation}

{\bf 7.} Now consider an arbitrary formula of the form (21). If the
numbers
$g_1$,\dots, $g_\epsilon$,
$s_1$,\dots, $s_\epsilon$,
$t_1$,\dots, $t_\epsilon$ satisfy  inequalities (20), then, as shown above,
formula (21) is equivalent to the formula
\begin{equation}
\newand_{\iota=1}^{\epsilon}
\forall y [-s_\iota^2-t_\iota^2-2<y\le s_\iota^2+t_\iota^2+2
\Rightarrow Z(g_\iota,s_\iota,t_\iota,y)>0].
%\eqno(30)
\end{equation}
We introduce the following notation:
$$
F_\iota(s_1,\dots,s_\iota,t_1,\dots,t_\iota)
=\sum_{\mu=1}^{\iota}
(2s_\mu^2+2t_\mu^2+4)
\quad (\iota=0,\dots,\epsilon),
$$
\begin{multline*}
Z_\iota(g_\iota, s_1,\dots,s_\iota,t_1,\dots,t_\iota,y)=\\
Z(g_\iota,s_\iota,t_\iota, %\\
y-F_{\iota-1} (s_1,\dots,s_{\iota-1},t_1,\dots,t_{\iota-1})-
s_\iota^2-t_\iota^2-2) \quad (\iota=1,\dots,\epsilon).
\end{multline*}
Obviously, formula (30) is equivalent to the formula
\begin{multline}
\newand_{\iota=1}^\epsilon
\forall y [
F_{\iota-1}(s_1,\dots,s_{\iota-1},t_1,\dots,t_{\iota-1})<y\le \\
\le F_\iota(s_1,\dots,s_\iota,t_1,\dots,t_\iota) \Rightarrow \\
\Rightarrow
Z_\iota(g_\iota, s_1,\dots,s_\iota,t_1,\dots,t_\iota,y)>0].
%\eqno(31)
\end{multline}

Let us denote by
$W_\epsilon(
g_1\dots,g_\epsilon, s_1,\dots,s_\epsilon,t_1,\dots,t_\epsilon,y)$
the polynomial
$$
\prod_{\iota=1}^\epsilon
Z_\iota(g_\iota, s_1,\dots,s_\iota,t_1,\dots,t_\iota,y).
$$

\

\noindent{\bf Lemma 4.}
If the numbers
$g_1\dots,g_\iota, s_1,\dots,s_\iota,t_1,\dots,t_\iota$
satisfy inequalities (20), then formula (21) is equivalent to the formula
\begin{multline*}
\forall f [
f\le F_{\epsilon}(s_1,\dots,s_{\epsilon},t_1,\dots,t_\epsilon)
\Rightarrow\\
\Rightarrow
W_\epsilon(g_1\dots,g_\epsilon, s_1,\dots,s_\epsilon,t_1,\dots,t_\epsilon,f)>0].
\end{multline*}

\

One can easily carry out the proof of the lemma using property (29).

\

\noindent{\bf Theorem 2.} Every recursively enumerable set of non-negative
integers can be represented by a formula of the form (2).

\

\noindent{\bf Proof.}
We will transform the formula
$
\exists z
{\mathfrak F}_2\,\&\,
{\mathfrak F}_3\,\&\,
{\mathfrak F}_4\,\&\,
{\mathfrak F}_5
$
into a form analogous to (21).

In the formula
${\mathfrak F}_2$ it suffices to transpose terms, which do
or do not contain $z$, to opposite sides of the inequalities.
We denote the resulting formula by ${\mathfrak F}_7$.

In the formula ${\mathfrak F}_3$
we replace the variable $d$, whose admissible values are nonnegative integers,
by the variable $z$, whose admissible values are all integers.
Since
$$
b\ge0,\quad (K(a,c))^{\lambda+1}>0,
$$
the formula thus obtained is equivalent to the formula
${\mathfrak F}_3$. Rewriting the formula we obtained
in a  form analogous to (18), we denote the new formula by
${\mathfrak F}_8$.

Analogously, in each conjunct of the formula ${\mathfrak F}_4$
we replace the variable $d$ by $z$. Since always
$$
b\ge0,\quad (K(a,c))^{(\lambda+1)^{\iota+1}}>
(c+1) (K(a,c))^{(\lambda+1)^{\iota}}>0,
$$
the resulting formula  is equivalent to the formula
${\mathfrak F}_4$.
We now perform  the same transformations on the formula thus obtained
as we had carried out with respect to the formula
${\mathfrak F}_4$ in the proof of Theorem~1. As a result,
we obtain a formula
${\mathfrak F}_9$, which is analogous to formula (19).

We replace the formula ${\mathfrak F}_5$ by an equivalent formula
\begin{multline*}
\exists z [ b-(c+1)(K(a,c))^{(\lambda+1)^\delta}
<
2(c+1)(K(a,c))^{(\lambda+1)^\delta}z< \\
(c+1)(K(a,c))^{(\lambda+1)^\delta}-b],
\end{multline*}
which we denote by ${\mathfrak F}_{10}$.

The formula
\begin{equation}
\exists z
{\mathfrak F}_7 \,\&\,
{\mathfrak F}_8 \,\&\,
{\mathfrak F}_9 \,\&\,
{\mathfrak F}_{10}
\end{equation}
is of a form analogous to (21).
The only difference is as follows:
the variables $g_\iota$, $s_\iota$, $t_\iota$ were replaced in (32)
by polynomials in the parameters $a$, $b$, $c$.
It is easy to check that for all values of the parameters,
the inequalities analogous to (20) hold. By Lemma 4 this enables us to find
the desired polynomials $F$ and $W$.

The theorem is proved.

\end{document}